\documentclass[12pt]{amsart}
\usepackage{amsmath,amssymb,amsthm}
\usepackage[]{latexsym}

\newcommand{\PP}{{\mathbb P}}

\newcommand{\cO}{{\mathcal O}}

\newtheorem{mtheorem}{Theorem}
\newtheorem{mlemma}[mtheorem]{Lemma}
\newtheorem{mproposition}[mtheorem]{Proposition}

\newcommand{\TT}{{\mathbb T}}

\newcommand{\Gal}{{\rm Gal}}

\newcommand{\sk}{\vspace{0.1in}}

\newcommand{\Aut}{\mbox{Aut}}
\newcommand{\PGL}{{\rm \bold{PGL}}}
\newcommand{\GL}{{\rm \bold{GL}}}
\newcommand{\End}{{\rm{End}}}

\newcommand{\SL}{{\rm {\bold SL}}}

\newcommand{\Z}{{\mathbb Z}}
\newcommand{\Q}{{\mathbb Q}}
\newcommand{\C}{{\mathbb C}}
\newcommand{\R}{{\mathbb R}}
\newcommand{\F}{{\mathbb F}}

\include{thebibliography}

\begin{document}

\title{On the group of automorphisms of \\ Shimura curves and applications\\
(Compositio Math. {\bfseries 131} (2002) 1-13)}

\author{Victor Rotger}

\footnote{Partially supported by a grant FPI from Ministerio de Educaci\'{o}n y
Ciencia, by DGCYT PB97-0893 and DGCYT PB96-0166.
}

\address{Dpt.\,Algebra i Geometria, Universitat de Barcelona,
Gran Via, 585, E-08007, Barcelona.}

\email{vrotger@mat.ub.es}



\begin{abstract}
Let $V_D$ be the Shimura curve over $\Q $ attached to the indefinite
rational quaternion algebra of discriminant $D$. In this note we
investigate the group of automorphisms of $V_D$ and prove that, in many
cases, it is the Atkin-Lehner group. Moreover, we determine the
family of bielliptic Shimura curves $(\mbox{over }\bar \Q \mbox{ and over }\Q
)$ and we use it to study the set of rational points on $V_D$ over quadratic
fields. Finally, we obtain explicit equations of
elliptic Atkin-Lehner quotients of $V_D$.
\end{abstract}

\maketitle

\section*{Introduction}
\sk\noindent

Let $B$ be an indefinite rational quaternion algebra and choose a maximal order
$\mathcal O \subset B$ of integers. Let $D=p_1\cdot ... \cdot p_{2r}$, $p_i$ prime numbers,
be the discriminant of $B$. We can view
$\Gamma =\{ \gamma \in \cO , n(\gamma )=1\} $ as an arithmetic subgroup of
$\SL _2(\R )$ through an identification $\Psi :B\otimes \R \cong M_2(\R )$ and
consider the Riemann surface $\Gamma \backslash \mathcal H $, where $\mathcal H $
denotes the upper-half plane of Poincar\'e. Shimura (\cite{Sh1}) showed that this is
the set of complex points of an algebraic curve $V_D/\Q $ over $\Q $ which
parametrizes abelian surfaces with quaternionic multiplication by $\cO $.

The classical modular case arises when we consider the split quaternion algebra
$B=M_2(\Q )$ of discriminant $D=1$. In this case, $V_1={\bf A} _{\Q }^1$ is the
$j$-line that classifies elliptic curves or, by squaring, abelian surfaces with
multiplication by $M_2(\Z )$. Throughout, we will limit ourselves to {\em non split}
quaternion algebras, that is, $D\not =1$. In this case, $\Gamma $ has no parabolic
elements and  $\Gamma \backslash \mathcal H $ is already compact so there are no
cusps and the automorphic
forms on $V_D$ do not admit Fourier expansions. In this regard, see \cite{Mor}.

As in the modular case, the elements of the Atkin-Lehner group
$W=\{ w_m ; m\mid D \} \cong C_2^{2r}$, where $C_2$ is the cyclic group of order
two, act as rational involutions on the Shimura curve $V_D$ and there
is a natural inclusion $W\subseteq \Aut _{\Q }(V_D)$ (see e.g. \cite{Jo2},
\cite{Ogg3}). The aim of the present
study is to examine the full group of automorphisms $\Aut (V_D\otimes\C )$ of these
curves. In the first section we describe its group structure and the field of
definition of its elements and we prove that $\Aut (V_D\otimes \C )=W$ in many
cases. In \cite{Ogg2}, Ogg studied the group of automorphisms of the modular
curves $X_0(N)$ for square-free levels $N$ and there, the action of
$\mbox{Aut}(X_0(N))$ on the set of cusps played a fundamental role.
When $D>1$, the difficulty lies precisely in the absence of cusps on $V_D$.

In section 2, the family of Shimura curves $V_D$ that admit bielliptic involutions
is determined. The hyperelliptic problem was already settled by Michon and Ogg
independently in \cite{Mi}, \cite{Ogg1} and \cite{Ogg3} and the family of bielliptic
modular curves $X_0(N)$ was given in \cite{Ba}. Since $\mbox{Aut}(X_0(N))$ is
largely understood (\cite{Ogg2}, \cite{KeMo}), the main point in \cite{Ba} was to
count the number of fixed points of the non
Atkin-Lehner involutions that appear when $4\mid N$ or $9\mid N$. In our case this
difficulty does not arise, but on the other hand the automorphism groups of the
Shimura curves $V_D$ are much less known.

In the last section, we derive some arithmetical consequences from the above results
concerning the set of rational points on $V_D$ over quadratic fields. Recall that
by a fundamental theorem of Shimura, there are no real points on Shimura curves and
therefore quadratic imaginary fields are the simplest fields over which these
curves may have rational points.  Our main theorem in this section
completely solves a question posed and studied by Kamienny in \cite{Ka}: which Shimura
curves $V_D$ of
genus $g\geq 2$ admit infinitely many quadratic points? This question is motivated
by Faltings' theorem on Mordell's conjecture and the answer is based upon ideas of
Abramovich, Harris and Silverman (see \cite{AbHa} and \cite{HaSi}).

Finally, we use the Cerednik-Drinfeld theory to compute equations
of elliptic Atkin-Lehner quotients of Shimura curves. Table 3 in
section 3 gives a Weierstrass equation of {\em all} elliptic
curves of the form $V_D/\langle w\rangle $ where $w\in \Aut(V_D)$
is any $\Q $-bielliptic involution on the curve. Some examples
were already given in \cite{Ro}.

The main tools used in this paper come from the reduction of
Shimura curves at both good and bad places. Drinfeld constructed a
projective model $M_D$ over $\Z $ of the Shimura curve $V_D$ which
extends the moduli interpretation given by Shimura to abelian
schemes over arbitrary bases (\cite{Dr}, \cite{Bu}, \cite{BoCa}).
Morita showed that $M_D$ has good reduction at any prime $p\nmid
D$ and Shimura (\cite{Sh1}) determined the zeta function of the
special fibre of $M_D$ at $p$. Moreover, the Cerednik-Drinfeld
theory (\cite{BoCa}, \cite{Ce}, \cite{Dr}, \cite{JoLi1}) provides
a good account of the behaviour of the reduction of $M_D \, $ (mod
$p$) when $p\mid D$.

{\em Acknowledgements. } I am indebted to Prof.\,Pilar Bayer for her
help and encouragement throughout this study. I also express my gratitude to
Prof.\,Gerald E.\,Welters for some helpful conversations, to Prof.\,Jordi Quer for
providing me with his algorithm on the decomposition of the Jacobian varieties of
modular curves and to the referee for the useful remarks. Finally, I thank the Tata
Institute of Fundamental Research and Mehta Research Institute for their warm
hospitality during the fall of 1999.

{\bfseries Notations.} We will denote by $\beta \mapsto \bar \beta $ the conjugation
map on $B$. The reduced trace and reduced norm on $B$ will be denoted respectively by
$\mbox{tr}(\beta )=\beta +\bar \beta $ and $\mbox{n}(\beta )=\beta \cdot \bar \beta $.

\section{The group of automorphisms of Shimura curves.}
\sk\noindent

Throughout, $V_D$ will denote the canonical model over $\Q $ of the Shimura curve of
discriminant $D=p_1\cdot ...\cdot p_{2r}\not =1$ (cf. \cite {Sh1}). It is a proper
smooth scheme over $\Q $ of dimension 1. Let $\Aut _{\Q }(V_D)$ be the group of
$\Q $-automorphisms of $V_D$ that sits inside the full group of automorphisms
$\Aut _{\C }(V_D \otimes \C )$ of the complex algebraic curve $V_D \otimes \C $.

\begin{mproposition}

If $g(V_D)\geq 2$, then
\begin{enumerate}
\item
All automorphisms of $V_D \otimes \C $ are defined over $\Q $.
That is: $\mbox{Aut} _{\C }(V_D\otimes \C )=\mbox{Aut} _{\Q }(V_D)$.

\item
$\mbox{Aut} _{\Q }(V_D)\cong C_2^s$, $s\geq 2r$.

\end{enumerate}
\end{mproposition}

{\em Proof:} In \cite{Ri1}, Ribet proved that all the endomorphisms of an abelian
variety $A/K$ with semistable reduction over a number field $K$ are defined over an
unramified extension of $K$. The Jacobian variety $J_D/\Q $ of $V_D$ has good
reduction at primes $p\nmid D$ and, from \cite{JoLi2}, we know that the identity's
connected
component of the reduction mod $p$, $p\mid D$, of the N\'eron model of $J_D$ is a
torus. Hence, $J_D$ has semistable reduction over $\Q $ and all its endomorphisms
are rational because $\Q $ has no non-trivial unramified extensions.

Since, by Hurwitz theorem, $\Aut _{\C }(V_D\otimes \C )$ is a finite group, any automorphism of
$V_D\otimes \C $ is defined over $\bar \Q $. Moreover, the natural map
$\Aut _{\bar \Q }(V_D\otimes \bar \Q )\rightarrow \Aut _{\bar \Q }
(J_D\otimes \bar \Q )$ is injective and $\Gal (\bar \Q/\Q )$-equivariant
and therefore we conclude from above that all automorphisms of $V_D\otimes \C $
are rational.

For the second part, let $X_0(D)/\Q $ be the modular curve of level $D$ and consider
the new part $J_0(D)^{new}/\Q $ of its Jacobian variety $J_0(D)$. It is well known
(\cite{Ri1})  that $\End^0_{\Q }(J_0(D)^{new})\cong \TT \otimes \Q \cong
\prod_{i=1}^tK_i$, where $\TT $ denotes the Hecke algebra of level $D$ and $K_i$ are
totally real number fields. Ribet's isogeny theorem (\cite{Ri2}) states the
existence of an isogeny

$$
\varphi : J_D \longrightarrow J_0(D)^{new}
$$
between $J_D$ and $J_0(D)^{new}$. This isogeny is Hecke invariant
(but sign-interchanging for the Atkin-Lehner action) and defined over $\Q $. Hence,
the ring of endomorphisms $\End _{\Q }(J_D) $ is an order in $\prod_{i=1}^t K_i$. An
automorphism of the curve $V_D$ induces an automorphism of finite order on $J_D$.
Moreover, the group of integral units in $\prod_{i=1}^t K_i$ is isomorphic to
$C_2^t$. So $\Aut _{\Q }(V_D)\cong C_2^s$ with $2r\leq s\leq t$, the first
inequality holding just because $W\subseteq \Aut _{\Q }(V_D)$.

\vspace{0.5cm}

We conclude that any automorphism of $V_D$ acts as a rational
involution on it. In view of the above proposition, we will simply denote  the group
$\Aut _{\C }(V_D\otimes \C )=\Aut _{\Q }(V_D)$ by $\Aut (V_D)$. Naturally we ask
whether the Atkin-Lehner group is the full group of automorphisms of the curve,
provided that $g(V_D)\geq 2$. This is the case for modular curves $X_0(N)$ of square
free level $N$, $N\not =37$ (\cite{Ogg2}, \cite{KeMo}) .

Recall that an elliptic point on the curve $V_D$ is a branched point of the natural projection $\mathcal H \rightarrow \Gamma \backslash \mathcal H$$=V_D(\C )$. The stabilizers of those elliptic points in $\Gamma/\{ \pm 1\} $  are of order 2 or 3. 2-elliptic points (resp. 3-elliptic points) correspond to $\Gamma $-conjugacy classes of embeddings of the quadratic order $\Z [i]$, $i^2=-1$ (resp. $\Z [\rho ]$, $\rho ^3=1$) in the quaternion order $\cO$. Their cardinality is given by

$$
e_2=\prod_{\ell \mid D} (1-{\underline {-4}\choose \ell }),
$$
$$
e_3=\prod_{\ell \mid D} (1-{\underline{-3} \choose \ell }),
$$
where $(\frac{\cdot }{\cdot })$ denotes the Kronecker symbol.

\begin{mtheorem}

Let $V_D$ be the Shimura curve of discriminant $D$. If it has no elliptic points,
then $\mathrm{Aut} (V_D)=W$.

\end{mtheorem}

{\em Proof:} If there are no elliptic points on $V_D(\C )$, then the natural
projection $\mathcal H \rightarrow \Gamma \backslash \mathcal H$$=V_D(\C )$ is the
universal cover of the Riemann surface $V_D(\C )$ so $\Aut (V_D)\cong
\mbox{Norm}_{\PGL _2^+(\R )}(\Gamma )/\Gamma $. Here the superindex $+$ denotes
matrices with positive determinant. It is known that
$W\cong \mbox{Norm}_{B^{\times }}(\Gamma )/\Q ^{\times }\Gamma \cong C_2^{2r}$
(\cite{Mi}, \cite{Ogg3}). We observe now that
the $\Q $-vector space spanned by $\Gamma $ is $\langle  \Gamma  \rangle _{\Q } =
B$. Indeed, since the reduced norm n is indefinite on the space of
{\em pure quaternions} $B_0=\{b\in B, \mbox{tr}(b)=0\} $, we can find linearly
independent elements $\omega _1$, $\omega _2$,
$\omega _3\in B_0$ such that $\Z [\omega _i]\subset B$ is a real quadratic order in
$B$. Then, by solving  the corresponding Pell equations, we find units $\gamma _i
\in \Z [\omega _i ]\cap \Gamma $, $\gamma _i \not =\pm 1$ such that $\{1, \gamma _1,
\gamma _2, \gamma _3\} $ is a $\Q $-basis of $B$.

Hence, any $\alpha \in \mbox{Norm}_{\GL _2^+(\R )}(\Gamma )$ will actually normalize
$B^{\times }$. By the Skolem-Noether theorem, $\alpha $ induces an inner
automorphism of $B$ so that $\alpha \in \R ^{\times}\mbox{Norm}_{B^{\times }}
(\Gamma )$. This shows that $\Aut (V_D)=W$.

{\em Remark: } In proving the above theorem, we have also
shown that the Atkin-Lehner group $W$ of an arbitrary Shimura curve
$V_D$ is exactly the subgroup of automorphisms that lift to a M\"{o}bius
transformation on $\mathcal H$ through the natural uniformization $\mathcal H
\rightarrow V_D(\C )=\Gamma \backslash \mathcal H$.

The proof remains valid for Eichler orders of square-free
level $N$ and therefore it generalizes an analogous result of Lehner and Newman for
discriminant $D=1$ (\cite{LeNe}).

The next theorem is similar in spirit to Theorem 2 and requires a previous lemma due to
Ogg.

\begin{mlemma}[\cite{Ogg2}]

Let $K$ be a field and $\mu (K)$ its group of roots of unity. Let
$p=\mathrm{max }(1, \mathrm{char }K)$ the characteristic exponent of $K$. Let $C$ be
an irreducible curve defined over $K$ and $P\in C(K)$ a regular point on it. Let $G$
be a finite group of $K$-automorphisms acting on $C$ and fixing the point $P$. Then
there is a homomorphism $f: G \rightarrow \mu (K)$ whose kernel is a $p$-group.

\end{mlemma}

\begin{mtheorem}

Let $D=2p$, $3p$; $p$ a prime number. If $g(V_D)\geq 2$, then $\mathrm{Aut} (V_D)=
W\cong C_2\times C_2$.

\end{mtheorem}

{\em Proof:} Suppose first that $D=2p$ with $p\equiv 3 \, $ (mod 4). In this case,
the fixed points on $V_D$ of the Atkin-Lehner involution $w_2$ are {\em Heegner
points} (see e.g. \cite{Al} for a
general account). Their coordinates on Shimura's canonical model $V_D$ generate
certain class fields. More precisely, if the genus $g(V_D)$ is even, then $w_2$
exactly fixes two points $P$, $P'$ with complex multiplication by the quadratic
order $\Z [i]$ and hence (\cite{ShTa}) $P$, $P'\in V_D(\Q (i))$. If $g(V_D)$ is odd,
then, besides $P$ and $P'$, $w_2$ fixes two more points $Q$,
$Q'\in V_D(\Q (\sqrt{ -2}))$  which have CM by $\Z [\sqrt{ -2}]$. As we have seen,
$\Aut (V_D)$ is an abelian group so it acts on the set of fixed points of $w_2$ on
$V_D$. Since all automorphisms are rational, they must keep the field of rationality
of these points so that $\Aut (V_D)$ actually acts on $\{ P, P'\} $. It follows from
the previous lemma that the order of the stabilizer of $P$ or $P'$ in $\Aut  (V_D)$
is at most 2. Hence, $\# \Aut (V_D)\leq 4$ and $\Aut (V_D)=W$.

Suppose now that $D=2p$, $p\equiv 1\, $ (mod 4) or $D=3p$,
$p\equiv 1\, $ (mod 3). By the theory of Cerednik and Drinfeld,
the special fibre $M_D \otimes \F _p$ of the reduction mod $p$ of
the integral model $M_D$ of our Shimura curve consists of two
rational irreducible components $Z$, $Z'$ defined over $\F
_{p^2}$. The complete local rings of the intersection points of
$Z$ and $Z'$ over the maximal unramified extension $\Z _p^{unr}$
of $\Z _p$ are isomorphic to $\Z _p^{unr}[x, y]/(x y - p^{\ell })$
for some {\em length} $\ell \geq 1$. The reduction mod $p$ of the
Atkin-Lehner involution $w_p$ switches $Z$ and $Z'$, fixing the
double points of intersection. Among these double points, there is
exactly one, say $\widetilde Q$, which has {\em length} 2, as it
follows from \cite{Ku}. Thus $\Aut (V_D)$ acting on $M_D\otimes \F
_p$ must fix $\widetilde Q$. Recalling now that $\Aut (V_D)\cong
C_2^s$, we again apply Ogg's lemma to the curve $Z/\F _{p^2}$
($p\not =2$) and the point $\widetilde Q$ to obtain that $\# \Aut
(V_D)\leq 4$. Therefore $\Aut (V_D)=W$.

In the remaining case, namely when $D=3p$, $p\equiv -1 \, $ (mod 3), we observe the
curious phenomenon that $\ell =109$ is a prime of good reduction for the Shimura
curve $V_D$ that yields

$$
\# M_D\otimes \F _{109} (\F _{109}) \not \equiv 0 \quad (\mbox{ mod } 4)
$$
except for the two exceptional cases $D=3\cdot 89$ and $D=3\cdot 137$. In any case,
we check that $\# M_{3\cdot 89}\otimes \F _{67} (\F _{67})=94$ and
$\# M_{3\cdot 137}\otimes \F _{103} (\F _{103})=98$. This is carried out by using
the explicit formula for the number of rational points over finite fields of the
reduction of Shimura curves at good places given by Jordan and Livn\'e in
\cite{JoLi1}. From this
we proceed as above: since all automorphisms of $V_D$ are defined over $\Q $, their
reduction mod $\ell $ must preserve the $\F _{\ell }$-rational points on
$M_D\otimes \F _{\ell }$ and we apply Ogg's lemma to the  regular curve
$M_D\otimes \F _{\ell }$ to conclude that $\Aut (V_D)=W$.

\vspace{0.5cm}
{\em Remark: } The first argument can be adapted for more general discriminants in an
obvious way. For instance, if $D=p\delta $ where $p$ is a prime integer,
$p\equiv 3 \, $ (mod 8), and ($\frac {-p}{\ell })=-1$ for any $\ell \mid \delta $,
then we again obtain that $\Aut (V_D)=W$ because the Hilbert class field of
$\Q (\sqrt{-p})$ is strictly contained in the ring class field of conductor 2 and,
by genus theory, $h(-p)$ is odd.

\vspace{0.5cm}
{\em Example: } Shimura curve quotient $V_{291}^+ = V_{291}/W$ has
genus 2 and therefore it is hyperelliptic. However, the hyperelliptic
involution on $V_{291}^+$ is exceptional: it does not lift to a M\"{o}bius
transformation on $\mathcal H$
through $\pi :\mathcal H \rightarrow
V_{291}^+(\C )=\Gamma \cdot W\backslash \mathcal H$, while all
automorphisms of $V_{291}$ are of Atkin-Lehner type by Theorem 4. This
is caused by the fact that $\pi $ is not the universal cover of
$V_{291}^+$.

\section{Bielliptic Shimura curves}
\sk\noindent

Recall that an algebraic curve $C$ of genus $g\geq 2$ is (geometrically)
bielliptic if it admits a degree 2 map $\varphi : C \rightarrow E$ onto a curve $E$
of genus 1. We will ignore fields of rationality until the next section.
Alternatively, $C$ is bielliptic iff there is an involutive automorphism acting on
it with $2g-2$ fixed points. We present now some facts about bielliptic curves $C$
such that, like Shimura curves, $\mathrm{Aut} (C)\cong C_2^s$.

\begin{mlemma} Let $C/K$, $\mbox{char }K \not =2$, be a bielliptic curve of genus
$g$ with $\mathrm{Aut} (C)\cong C_2^s$ for some $s\geq 1$. For any $w\in \mathrm{Aut} (C)$, let $n(w)$
denote the number of fixed points of $w$ on $C$. Let $v$ be a bielliptic involution
on $C$ and for any $w\in \mathrm{Aut} (C)$, $w\not =1 \mbox{ or } v$, denote
$w'=v\cdot w$.

\begin{enumerate}
\item
If $g$ is even, then $n(w)=2$ and $n(w')=6$, or vice versa. If $g$ is odd, then
$\{ n(w), n(w')\} =\{ 0, 0\}, \{0, 8\} $ or $\{ 4, 4\} $ as non-ordered pairs.

\item
If $g$ is even, then $s\leq 3$. If $g$ is odd, then $s\leq 4$.

\item
If $g\geq 6$, then the bielliptic involution $v$ is unique.

\end{enumerate}
\end{mlemma}

{\em Proof: } The first part follows from Hurwitz's theorem applied to the map
$C \rightarrow C/\langle v, w\rangle $, while $2.$ and $3.$ are simple corollaries
of that.

\vspace{0.5cm}
{\em Remark: }Observe that if $D$ is odd, then $g(V_D)$ is always odd, as we check
from Eichler's formula for the genus (see e.g. \cite{Ogg3}).

\vspace{0.5cm}
Obviously, the main source for possible bielliptic involutions on the curves $V_D$
is the Atkin-Lehner group. From Eichler's formula for $n(w)$, $w\in W$
(see \cite{Ogg3}), it
is a routine exercise to check whether $V_D$ has bielliptic involutions of
Atkin-Lehner type. An alternative way to compute $n(w)$ is to read {\em backwards}
the last column of Table 5 in \cite{Ant}. This is because Ribet's isogeny
$\varphi : J_D \rightarrow J_0(D)^{new}$ switches the sign of the Atkin-Lehner
action. But, first, we should focus on possible extra involutions and also bound the
bielliptic discriminants $D$. Following Ogg's method in \cite{Ogg1}, we give such an
upper bound in the next

\begin{mproposition}
If $D>547$, $V_D$ is not bielliptic.

\end{mproposition}

{\em Proof: } Suppose that the curve $V_D$ is bielliptic: there is a degree 2 map
$\varphi : V_D\rightarrow E$ onto a curve $E$ of genus 1. By Proposition 1.2, both
$\varphi $ and $E$ are defined over $\Q $ although $E$ may not be an elliptic curve
over $\Q $ since it may fail to have rational points (see Section 3 for examples).
Choose a prime of good reduction $\ell \nmid D$ of $V_D$, let $K_{\ell }$ be
the
quadratic unramified extension of $\Q _{\ell }$ and let $R_{\ell }$ denote its
ring of integers. As follows from \cite{JoLi1},
$V_D(K_{\ell })\not =\emptyset $ and hence $E$ is an elliptic curve over $K_{\ell
}$. Moreover, due to Ribet's isogeny theorem,
$E$ also has good reduction over $\ell $. By the universal property of the N\'eron
model of $E$ over $R_{\ell }$, $\varphi $ extends to the minimal smooth model
$M_D\otimes R_{\ell }$ of $V_D$  and we can reduce the bielliptic structure
mod $\ell $ to obtain a $2:1$
map $\tilde \varphi : M_D\otimes \F _{\ell ^2} \rightarrow \widetilde E$. From
Weil's estimate, $N_{\ell ^2}=\# M_D\otimes \F _{\ell ^2}(\F _{\ell ^2})\leq 2\cdot
\# \widetilde E(\F _{\ell ^2})\leq 2(\ell +1)^2$.  Besides, we obtain from
\cite{JoLi1} that $\frac {(\ell -1)}{12}\prod _{p\mid D}{(p-1)}\leq N_{\ell ^2}$
so $N_{\ell ^2}$ grows as $D$ tends to infinity. Applying these inequalities for
$\ell =2, 3, 5, 7, 11$, we conclude that, if $2\cdot3\cdot5\cdot7\cdot11\nmid D$,
then $D\leq 546$. But, if $2\cdot 3\cdot 5\cdot 7\cdot 11\mid D$,
then $s\geq 5$ and therefore $V_D$ cannot be bielliptic (Lemma 5.2).

\vspace{0.5cm}
We are now able to prove

\begin{mtheorem}

There are exactly thirty-two values of $D$ for which $V_D$ is bielliptic. In each case, the bielliptic involutions are of Atkin-Lehner type. These values, together with the genus $g=g(V_D)$ and the bielliptic involutions are given in Table 1 below.

\newpage
\vspace{1cm}
\begin{tabular}{|r|r|l|r|r|r|l|r|r|r|l|} \hline
\multicolumn{11}{|c|}{\bfseries Table 1: Bielliptic Shimura curves}\\

\hline
D&g&$w_m$&   &D&g&$w_m$&   &D&g&$w_m$ \\
\hline \hline
$26$&$2$&$w_2$,$w_{13}$&  &$82$&$3$&$w_{82}$&   &$210$&$5$&$w_{30}$,$w_{42}$, \\
$35$&$3$&$w_7$&   &$85$&$5$&$w_{17}$&   &   &   &$w_{70}$,$w_{105}$, \\
$38$&$2$&$w_2$,$w_{19}$&   &$94$&$3$&$w_2$&   &   &   &$w_{210}$ \\
$39$&$3$&$w_{13}$&   &$106$&$4$&$w_{53}$,$w_{106}$&  &$215$&$15$&$w_{215}$  \\
$51$&$3$&$w_3$&   &$115$&$6$&$w_{23}$&   &$314$&$14$&$w_{314}$   \\
$55$&$3$&$w_5$&   &$118$&$4$&$w_{59}$,$w_{118}$&   &$330$&$5$&$w_3$,$w_{22}$  \\
$57$&$3$&$w_{57}$&    &$122$&$6$&$w_{122}$&   &   &   &$w_{33}$,$w_{165}$, \\
$58$&$2$&$w_2$,$w_{58}$&    &$129$&$7$&$w_{129}$&    &    &    &$w_{330}$ \\
$62$&$3$&$w_2$&   &$143$&$12$&$w_{143}$&   &$390$&$9$&$w_{390}$ \\
$65$&$5$&$w_{65}$&   &$166$&$6$&$w_{166}$&   &$462$&$9$&$w_{154}$ \\
$69$&$3$&$w_3$&    &$178$&$7$&$w_{89}$&   &$510$&$9$&$w_{510}$ \\
$77$&$5$&$w_{11}$,$w_{77}$&   &$202$&$8$&$w_{101}$&   &$546$&$13$&$w_{546}$\\
\hline
\end{tabular}

\end{mtheorem}
\newpage

{\em Proof: } Since we need only consider discriminants $D\leq 546$, we can first
use any programming package to build up the list of Atkin-Lehner bielliptic
involutions on Shimura curves $V_D$. These computations yield Table $1$ above.
In order to ensure that no extra bielliptic involutions arise, we observe that the
above results (and in particular Theorem 4) imply that any bielliptic involution on
$V_D$, for most of the discriminants $D\leq 546$, must be of Atkin-Lehner type.
There are exactly three cases, namely $D=55$, $D=85$ and $D=145$, for which none of
the previous results and their obvious generalizations seem to apply.

{\em Ad hoc} arguments can be worked out for them. Firstly, the Jacobian varieties
of the curves $V_{55}$ and $V_{85}$ have just one $\Q $-isogeny class of sub-abelian
varieties of dimension 1, so there can be at most one bielliptic involution on these
curves. But $w_5$ (resp. $w_{17}$) is already a bielliptic involution on $V_{55}$
(resp. $V_{85}$).

More interesting is the curve $V_{145}$ of genus 9. It is not bielliptic by any
Atkin-Lehner involution although $J_{145}\sim _{\Q }E\times S\times A_3\times A'_3$,
where each factor has dimension 1, 2, 3 and 3 respectively. We check that
$n(w_5)=n(w_{29})=0$ and $n(w_{145})=8$, so if there was a bielliptic involution $v$
on $V_{145}$ then $n(w'_{145})=0$,  by Lemma 5.1. It follows from Lefschetz's fixed
point formula (see \cite{BiLa}) that the rational traces of these three involutions
on the Jacobian $J_{145}$ would be $\mbox{tr}(w_5)=\mbox{tr}(w_{29})=
\mbox{tr}(w'_{145})=2$.  Moreover, involutions on $J_{145}$ must be of the form
$\{ \pm 1_{E}\}\times \{ \pm 1_S\} \times \{ \pm 1_{A_3}\} \times
\{ \pm 1_{A'_3}\} $, up to conjugation by $\varphi $, thus $\mbox{tr}=2$ can only be
attained by two different involutions. Therefore $v$ cannot exist and $V_{145}$ is
not bielliptic.

It can be showed that actually $\Aut (V_{145})=W$: from the decomposition of
$J_{145}$ we know that $W\cong C_2^2\subseteq \Aut (V_{145})\subseteq C_2^4$. Since
$V_{145}$ is neither hyperelliptic (\cite{Ogg3}) nor bielliptic
(as we have just seen), it follows that the involutions
$\{ -1_{E}\}\times \{ -1_S\} \times \{ -1_{A_3}\} \times \{ -1_{A'_3}\} $ and
$\{ +1_{E}\}\times \{ -1_S\} \times \{ -1_{A_3}\} \times \{ -1_{A'_3}\} $ cannot be
induced from $\Aut (V_{145})$. Thus it is a subgroup of index at least 4 in $C_2^4$
and $\Aut (V_{145})=W$.

\section{Infinitely many quadratic points on Shimura curves}
\sk\noindent

In \cite{Sh2}, Shimura proved that  $V_D(\R )=\emptyset $ and in particular there
are no $\Q $-rational points on Shimura curves $V_D$. Jordan and Livn\'e
(\cite{JoLi1}) gave explicit criteria for deciding whether the curves $V_D$ do have
rational points over the $p$-adic fields $\Q _p$ for any finite prime $p$.

Much less is known about rational points over global fields. Jordan (\cite{Jo1})
proved that for a fixed quadratic imaginary field $K$, with class number
$h_K\not =1$, there are only finitely many discriminants $D$ for which $K$ splits
the quaternion algebra $B$ of discriminant $D$ and $V_D(K)\not =\emptyset $. In this
section we solve a question that is to an extent reciprocal: which Shimura curves
$V_D$, $g(V_D)\geq 2$, have infinitely many quadratic points over $\Q $?

That is,
$$
\# V_D(\Q , 2)=\# \{ P\in V_D(\bar \Q ), [ \Q (P):\Q ] \leq 2\}=+\infty
$$
We will say that an
algebraic curve $C/K$ of genus $g\geq 2$ is {\em hyperelliptic} over $K$
(respectively {\em bielliptic} over $K$) if there is an involution $v\in \Aut
_K(C)$ such that the quotient curve $C/\langle v\rangle
$ is $K$-isomorphic to $\PP _{K }^1$ (resp. an elliptic curve $E/K$).
Notice that in both cases $C/\langle v\rangle (K)\not =\emptyset $
while it perfectly well happen that $C(K)=\emptyset $.

The following theorem of Abramovich and Harris shows that the
question above is closely related to the diophantine
problem of determining the family of
hyperelliptic and bielliptic Shimura curves over $\Q $.

\begin{mtheorem}[\cite{AbHa}]

Let $C$ be an algebraic curve of genus greater than or equal to $2$, defined over a
number field $K$. Then $C(K , 2)=\# \infty $ if and only if $C$ is either
hyperelliptic over $K $ or bielliptic over $K $ mapping to an elliptic curve $E$
of $\mbox{rank} _{K }(E)\geq 1$.

\end{mtheorem}

Ogg (\cite{Ogg3}, \cite{Ogg4}) gave the list of hyperelliptic Shimura curves over
$\Q $. In what follows, we will determine which bielliptic Shimura curves from
Table 1 are bielliptic over $\Q $.

We first observe that the map $V_D \rightarrow V_D/\langle w\rangle $ is always
defined over $\Q $ since we showed that
$\Aut _{\C }(V_D\otimes \C )=\Aut _{\Q }(V_D)$ (Proposition 1). In order to
check whether $V_D/\langle w\rangle (\Q )\not =\emptyset $ for each pair $(D$,$w)$
in Table 1, we can disregard those in which $V_D/\langle w\rangle $ fails to have
rational points over some completion $\Q _v$ of $\Q $. This is done by using
the precise results in that direction given by Jordan and Livn\'e in  \cite{JoLi1}
and Ogg in \cite{Ogg3} and \cite{Ogg4}; the conclusions are compiled in the
following table. Let us say that a field $L$ is {\em deficient}
for an algebraic curve $C$
defined over a subfield $K\subset  L$ if $C(L)=\emptyset $.

\newpage
\vspace{1cm}
\begin{tabular}{|l|l|l|l|l|l|l|l|l|l|l|}
\hline
\multicolumn{11}{|c|}{\bfseries Table 2: Deficient completions $L$
of $\Q $ for $V_D/\langle w_m\rangle $}\\

\hline
D&$w_m$&L&&D&$w_m$&L&&D&$w_m$&L \\
\hline \hline

$35$&$w_7$&$\Q _5$& &$115$&$w_{23}$&$\Q _5$& &$330$&$w_3$&$\R $,$\Q _2$\\

$39$&$w_{13}$&$\R $,$\Q _3$& &$178$&$w_{89}$&$\R $,$\Q _2$& &  &  &$\Q _5$,$\Q _{11}$\\

$51$&$w_{3}$&$\Q _{17}$& &$210$&$w_{30}$&$\R $,$\Q _3$& &$330$&$w_{22}$&$\R $,$\Q _2$,$\Q _3$\\

$55$&$w_5$&$\R $,$\Q _{11}$& &$210$&$w_{42}$&$\R $,$\Q _2$,$\Q _3$& &  &  &$\Q _5$,$\Q _{11}$\\

$62$&$w_2$&$\R $,$\Q _{31}$& & & &$\Q _5$,$\Q _7$& &$330$&$w_{33}$&$\R $,$\Q _2$\\

$69$&$w_3$&$\R $,$\Q _{23}$& &$210$&$w_{70}$&$\R $,$\Q _2$&  &&&$\Q _3$,$\Q _5$\\

$77$&$w_{11}$&$\R $,$\Q _7$& &  & &$\Q _3$,$\Q _5$ &  &$330$&$w_{165}$&$\Q _2$,$\Q _3$\\

$85$&$w_{17}$&$\Q _5$& &$210$&$w_{105}$&$\R $,$\Q _2$& & & &$\Q _5$,$\Q _{11}$\\

$94$&$w_2$&$\R $,$\Q _{47}$& & & &$\Q _7$& &$462$&$w_{154}$&$\R $,$\Q _{11}$\\

\hline
\end{tabular}
\newpage

On the genus 1 Atkin-Lehner quotients $V_D/\langle w_m\rangle $  that do have
rational points over all completions of $\Q $, we can try to construct a
$\Q $-rational point by means of the theory of complex multiplication.
That is, a Heegner point $P\in V_D(K)$ with CM by a quadratic imaginary order $R$,
$R\otimes \Q =K$, $h(R)=1$, will project onto a $\Q $-rational point on
$V_D/\langle w_m\rangle $ if and only if $w_m(P)=\bar P$, where $\bar P$ is the
complex
conjugate of $P$ on $V_D(K)$. From \cite{Jo2}, 3.1.4, we deduce that $w_m(P)=\bar P$
if $m$ is the product of the primes $p\mid D$ that are inert in $K$.

Performing the necessary computations, it follows that among those pairs $(D, w)$
that $V_D/\langle w\rangle (\Q _v)\not =\emptyset $ for every completion $\Q _v$ of
$\Q $, it is always possible to produce a $\Q $-rational point on
$V_D/\langle w\rangle $ by the above means, except for two interesting cases:
$(V_{26}$, $w_2)$ and $(V_{58}$, $w_2)$.

Since $\mbox{g}(V_{26})=\mbox{g}(V_{58})=2$, we may apply a result of
Kuhn (\cite{Kuhn}) to deduce that there are also rational points on the
quotients $V_{26}/\langle w_2\rangle $ and $V_{58}/\langle w_{2}\rangle
$. Therefore, the Hasse-Minkowsky principle is never violated for the
Atkin-Lehner quotients from Table 1 and those pairs $V_D/\langle w_m\rangle
$ that do not appear in Table 2 are bielliptic curves over $\Q $. There
are only eighteen values of $D$ for which $V_D$ admits a bielliptic involution over
$\Q $.

It still remains to compute the Mordell-Weil rank of the elliptic curves
$V_D/\langle w\rangle $ over $\Q $. Using Cremona's tables \cite{Cre}, switching the
sign of the Atkin-Lehner action as explained above, we can determine the
$\Q $-isogeny class of these elliptic curves.

This is enough to compute their Mordell-Weil rank but we can use a
beautiful idea of Roberts (\cite{Ro}) to compute the $\Q
$-isomorphism class and hence a Weierstrass equation for them as
follows: Cremona's tables give the Kodaira symbols of the
reduction of elliptic curves $E$ at the primes $p\mid \mbox{cond}
(E)$. This is done by using Tate's algorithm which makes use of a
Weierstrass equation of the curve. This is not available in our
case, but we can instead use the Cerednik-Drinfeld theory to
compute the Kodaira symbols for the reduction mod $p$, $p\mid D$,
of $V_D/\langle w\rangle $ and contrast them with Cremona's
tables. This procedure uniquely determines the $\Q $-isomorphism
class of the curves.

\vspace{0.5cm}
{\em Example: } Curve $V_{210}$ has genus 5 and is bielliptic by the Atkin-Lehner
involution $w_{210}$. After Eichler's theory on optimal embeddings (see e.g.
\cite{Ogg3}), the quadratic order $\Z [\sqrt{-43} ]$ embeds in the quaternion
algebra $B$ of discriminant 210. Such an embedding produces a point
$P\in V_{210}(\Q (\sqrt{-43} ))$. From the above, it follows that
$w_{210}(P)=\bar P$. Therefore, $V_{210}/\langle w_{210}\rangle (\Q )\not
=\emptyset $ and we obtain that $(V_{210}, w_{210})$ is a bielliptic pair over
$\Q $. A glance at Cremona's Table 3, pp. 249-250, shows that the elliptic curve
$V_{210}/\langle w_{210}\rangle $ falls in the $\Q $-isogeny class $210D$ because it
is the only one that corresponds to a newform $f\in H^0(\Omega ^1,J_{210})$ such
that $w_{210}^*(f)=f$ (recall that the sign for the Atkin-Lehner action is opposite
to the classical modular case!). Therefore, from Cremona's Table 4,
$\mbox{rank} _{\Q }(V_{210}/\langle w_{210}\rangle )=1$.

In order to determine a Weierstrass equation for $V_{210}/\langle
w_{210}\rangle $ we may compute the Kodaira symbols of its
reduction mod $p$, $p\mid 210$. It suffices to study the reduction
at $p=3$. The Cerednik-Drinfeld theory asserts that
$M_{210}\otimes \F _3$ is reduced and its irreducible components
are all rational and defined over $\F _9$. Moreover,
$M_{210}\otimes \Z _3$ is a (minimal) regular model over $\Z _3$.
This is because over the quadratic unramified integral extension
$R_3$ of $\Z _3$, $M_{210}\otimes R_3$ is a Mumford curve
uniformized by a (torsion-free) Schottky group, as one checks from
Cerednik-Drinfeld's explicit description of this group and the
congruences $5\equiv -1\quad (\mbox{mod } 3)$ and $7\equiv -1\quad
(\mbox{mod } 4)$.

In a way, Cerednik-Drinfeld's description of the reduction of
Shimura curves at $p\mid D$ is not so different from
Deligne-Rappoport's for the modular curves $X_0(N)$ at $p\parallel
N$ because $M_{p\delta }\otimes \F _p$ is again the union of two
copies of the Shimura curve -also called Gross curve- $M_{\delta
}\otimes \F _p$, defined in terms of a {\em definite} quaternion
algebra.

Let $h(\delta ,\nu )$ denote the class number of an (arbitrary) Eichler order of
level $\nu $ in the quaternion algebra of discriminant $\delta $. The dual graph
$\mathcal G$ of $M_{210}\otimes \F _3$ has as vertices the irreducible components of
$M_{210}\otimes \F _3$. There are $2h(\frac{210}{3}, 1)=2h(70, 1)=4$ of them. Two
vertices $v$, $\tilde v$ in $\mathcal G$ are joined by as many edges as there are
intersection points between the corresponding components $Z$, $\widetilde Z$ in
$M_{210}\otimes \F _3$. In our case, there are $h(\frac{210}{3}, 3)=8$ edges in
$\mathcal G$, that is, 8 double points in $M_{210}\otimes \F _3$.

We may label the 4 vertices $v_1$, $v_1'$, $v_2$, $v_2'$ so that $w_3(v_i)=v_i'$,
where $w_3$ still denotes the involution $w_3$ now acting on $\mathcal G$. There are
no edges joining $v_1$ and $v_2$, and the same holds for $v_1'$ and $v_2'$. The
total number of edges joining $v_1$ with $v_1'$ and $v_2$ with $v_2'$ is 4, as
Kurihara (\cite{Ku}) deduced from trace formulae of Brandt matrices. Since there
must also be $p+1=4$ edges at the star of any vertex, it turns out that the dual
graph $\mathcal G$ must be

\quad
\quad
\quad
\quad
\quad
\quad
\quad

\begin{picture}(350,80)

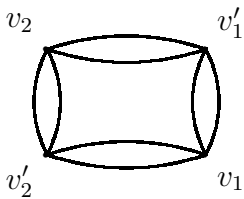
\put(300,10){\makebox(0,0){\bfseries Dual graph of $M_{210}\otimes \F _3$  }}
\put(120,20){\circle*{2}}
\put(180,20){\circle*{2}}
\put(120,60){\circle*{2}}
\put(180,60){\circle*{2}}

\put(110,10){\makebox(0,0){$v_2'$}}
\put(190,10){\makebox(0,0){$v_1$}}
\put(110,70){\makebox(0,0){$v_2$}}
\put(190,70){\makebox(0,0){$v_1'$}}

\thicklines\qbezier(120,20)(110,40)(120,60)
\qbezier(120,20)(130,40)(120,60)
\qbezier(120,20)(150,10)(180,20)
\qbezier(120,20)(150,30)(180,20)
\qbezier(180,60)(150,70)(120,60)
\qbezier(180,60)(150,50)(120,60)
\qbezier(180,60)(170,40)(180,20)
\qbezier(180,60)(190,40)(180,20)

\end{picture}

Since $3\mid 210$, $w_{210}(\{ v_1, v_2\} )=\{v_1', v_2'\} $ and therefore
$\mathcal G$$/\langle w_{210}\rangle $ is a graph with two vertices joined by two
edges, which corresponds to the Kodaira symbol $I_2$. The only elliptic curve in the
$\Q $-isogeny class $210D$ whose reduction type at $p=3$ is $I_2$ is $210D2$. Hence,
a Weierstrass equation for $V_{210}/\langle w_{210}\rangle $ is
$y^2+xy=x^3+x^2-23x+33$.

Performing similar computations, we obtain the list of bielliptic Shimura curves
$(V_D$, $w)$ over $\Q $ such that the genus 1 Atkin-Lehner quotient
$V_D/\langle w\rangle $ is an elliptic curve with infinitely many rational points.
With this procedure, we also give a Weierstrass equation for the elliptic curves
$V_D/\langle w\rangle $. Together with the hyperelliptic Shimura curves over $\Q $
given by Ogg, we obtain the family of Shimura curves of genus $g(V_D)\geq 2$ with
infinitely many quadratic points. Summing up, we obtain the following

\begin{mtheorem}
There are only finitely many $D$ for which $V_D$ has infinitely many quadratic points over $\Q $. These curves, together with their rational or elliptic quotients, are listed below.
\newpage
\vspace{1cm}
\begin{tabular}{|c|c|c|c|c|c|c|c|c|c|c|}
\hline
\multicolumn{11}{|c|}{\bfseries Table 3: Shimura curves $V_D$, $g(V_D)\geq 2$, with $\# V_D(\Q, 2)=+\infty $}\\
\hline
D&$w_m$&$V_D/\langle w_m\rangle $&   &D&$w_m$&$V_D/\langle w_m\rangle $&   &D&$w_m$&$V_D/\langle w_m\rangle $ \\
\hline \hline
$26$&$w_{13}$& $\PP _{\Q }^1$ &  &$77$&$w_{77}$& $77A_1$ &   &$143$&$w_{143}$ &$143A_1$ \\
$35$&$w_{35}$ &$\PP _{\Q }^1$&   &$82$&$w_{82}$& $82A_1$ &   &$146$&$w_{146}$&$\PP _{\Q }^1$ \\
$38$&$w_{38}$ &$\PP _{\Q }^1$&   &$86$ & $w_{86}$ &$\PP _{\Q }^1$ &   &$159$&$w_{159}$&$\PP _{\Q }^1$ \\
$39$&$w_{39}$ &$\PP _{\Q }^1$ &   &$87$&$w_{87}$ & $\PP _{\Q }^1$&  &$166$&$w_{166}$&$166A_1$  \\
$51$&$w_{51}$ &$\PP _{\Q }^1$&   &$94$&$w_{94}$& $\PP _{\Q }^1$&   &$194$&$w_{194}$&$\PP _{\Q }^1$   \\
$55$&$w_{55}$ & $\PP _{\Q }^1$&   &$95$&$w_{95}$ &$\PP _{\Q }^1$ &   &$206$&$w_{206}$ &$\PP _{\Q }^1$  \\
$57$&$w_{57}$ &$57A_1$&     &$106$&$w_{106}$& $106B_1$&   &$210$&$w_{210}$ &$210D_2$ \\
$58$&$w_{29}$& $\PP _{\Q }^1$&    &$111$&$w_{111}$&$\PP _{\Q }^1$&    &$215$&$w_{215}$&$215A_1$ \\
   &$w_{58}$&$58A_1$&   &$118$&$w_{118}$ & $118A_1$ &   &$314$&$w_{314}$&$314A_1$ \\
$62$&$w_{62}$ & $\PP _{\Q }^1$&    &$119$&$w_{119}$ & $\PP _{\Q }^1$&   &$330$&$w_{330}$ & $330E_2$ \\
$65$&$w_{65}$ &$65A_1$&    &$122$&$w_{122}$ & $122A_1$&   &$390$&$w_{390}$& $390A_2$ \\
$69$&$w_{69}$ & $\PP _{\Q }^1$&    &$129$&$w_{129}$ & $129A_1$&    &$510$&$w_{510}$& $510D_2$ \\
$74$&$w_{74}$ & $\PP _{\Q }^1$&   &$134$&$w_{134}$&$\PP _{\Q }^1$ &   &$546$&$w_{546}$& $546C_2$   \\
\hline
\end{tabular}

\end{mtheorem}
\pagebreak

\end{document}